\documentclass[a4paper,draft]{amsart}
\usepackage{amsmath}
\usepackage{amssymb}
\usepackage{bbm}
\usepackage{a4wide}

\parskip\smallskipamount

\def\<#1>{\langle#1\rangle}

\def\Nice{\mathrm{Nice}}
\def\qbinom#1#2{\genfrac{[}{]}{0pt}{}{#1}{#2}}
  \font\tworm=cmr3

\newcommand\todo[1][.]{\edef\tmpa{.}\edef\tmpb{#1}%
  \ifx\tmpa\tmpb
    \typeout{To Be on page \thepage}\fbox{\bf To Be}
  \else
    \typeout{To Be on page \thepage: #1}\fbox{{\bf To Be:} #1}
  \fi
}

\begin{document}

 \author[Manuel Kauers]{Manuel Kauers\,$^\ast$}
 \address{Manuel Kauers, Research Institute for Symbolic Computation, J. Kepler University Linz, Austria}
 \email{mkauers@risc.uni-linz.ac.at}
 \thanks{$^\ast$ Supported in part by the Austrian Science Foundation (FWF)
  grants P19462-N18 and P20162-N18.}

 \author[Christoph Koutschan]{Christoph Koutschan\,$^\ast$}
 \address{Christoph Koutschan, Research Institute for Symbolic Computation, J. Kepler University Linz, 
        Austria}
 \email{ckoutsch@risc.uni-linz.ac.at}

 \author[Doron Zeilberger]{Doron Zeilberger\,$^{\ast\ast}$}
 \address{Doron Zeilberger, Mathematics Department, Rutgers University (New Brunswick), Piscataway, NJ, USA.}
 \email{zeilberg@math.rutgers.edu}

 \thanks{$^{\ast\ast}$ Supported in part by the United States of America National Science Foundation.}

 \title[A Proof of George Andrews' and Dave Robbins' $q$-TSPP Conjecture]
{A~Proof of George~Andrews' and Dave~Robbins' $q$-TSPP~Conjecture {\tworm(modulo a finite amount of routine calculations)}}

 \maketitle

{\small
Accompanied by Maple packages TSPP and qTSPP available from
\hfill\break
\verb|http://www.math.rutgers.edu/~zeilberg/mamarim/mamarimhtml/qtspp.html|.}

\bigskip
\qquad\qquad\qquad\qquad\qquad\qquad
\qquad\qquad\qquad\qquad\qquad\qquad {\it Pour Pierre Leroux, In Memoriam}
\bigskip

\section*{Preface: Montr\'eal, May 1985}

In the historic conference {\it Combinatoire \'Enum\'erative}~\cite{LL}
wonderfully organized by Gilbert Labelle and {\bf Pierre Leroux}
there were many stimulating lectures, including a very interesting one
by Pierre Leroux himself, who talked about his joint work with
Xavier Viennot~\cite{LV}, on solving differential equations
combinatorially! 
During the problem session of that very same \emph{colloque},
chaired by Pierre Leroux,
Richard Stanley raised some intriguing problems about the
enumeration of plane partitions, that he later expanded
into a fascinating article~\cite{Sta1}. Most of these problems
concerned the enumeration of \emph{symmetry classes} of
\emph{plane partitions}, that were discussed in more detail
in another article of Stanley~\cite{Sta2}.  \emph{All} of the
conjectures in the latter article have since been proved
(see Dave Bressoud's modern classic~\cite{B}), \emph{except}
one, that, so far, \emph{resisted} the efforts of the greatest minds
in enumerative combinatorics. It concerns the proof of an
explicit formula for the $q$-enumeration
of \emph{totally symmetric plane partitions}, conjectured independently by
George Andrews and Dave Robbins~(\cite{Sta2}, \cite{Sta1} (conj. 7), \cite{B} (conj.~13)). 
In this tribute to Pierre Leroux, we describe how to prove that last stronghold.

\section{$q$-TSPP: The Last Surviving Conjecture About Plane Partitions}

Recall that a \emph{plane partition} $\pi$ is an array
$\pi=(\pi_{ij})$, ${i,j \geq 1}$, of positive integers
$\pi_{ij}$ with finite sum $|\pi|=\sum \pi_{ij}$,
which is weakly decreasing in rows and columns.
By stacking $\pi_{ij}$ unit cubes on top of the
$ij$ location, one gets the 3D Ferrers diagram,
that can be identified with the plane-partition,
and is a left-, up-, and bottom- justified structure
of unit cubes, and we can refer to the locations
$(i,j,k)$ of the individual unit cubes.

A plane partition is \emph{totally symmetric} iff
whenever $(i,j,k)$ is occupied (i.e. $\pi_{ij} \geq k$),
it follows that all its (up to $5$) permutations:
$\{(i,k,j),(j,i,k),(j,k,i),(k,i,j),(k,j,i)\}$ are
also occupied. In 1995, John Stembridge~\cite{Ste} proved 
Ian Macdonald's conjecture that
the number of totally symmetric plane partitions (TSPPs)
whose $3D$ Ferrers diagram is bounded inside the cube
$[0,n]^3$ is given by the nice product-formula
$$
\prod_{ 1 \leq i \leq j \leq k \leq n} \frac{i+j+k-1}{i+j+k-2}.
$$
Ten years after Stembridge's completely human-generated proof,
George Andrews, Peter Paule, and Carsten Schneider~\cite{APS} came up
with a \emph{computer-assisted} proof, that, however required lots
of human ingenuity and ad hoc tricks, in addition to a considerable
amount of computer time. 

Way back in the early-to-mid eighties (ca. 1983),
George Andrews and Dave Robbins independently conjectured
a \emph{$q$-analog} of this formula, namely that the
\emph{orbit-counting generating function} (\cite{B}, p.~200, \cite{Sta1}, p.~289)
is given by
$$
\prod_{1 \leq i \leq j \leq k \leq n} 
\frac{1-q^{i+j+k-1}}{1-q^{i+j+k-2}}.
$$

In this article we will show how to prove this conjecture
(modulo a finite amount of routine computer calculations that
may be already feasible today [with great technical effort],
but that would most likely be routinely checkable on a standard
desktop in twenty years).

\section{Soichi Okada's Crucial Insight}

Our starting point is an elegant reduction, by Soichi Okada~\cite{O},
of the $q$-TSPP statement, to the problem of evaluating
a certain ``innocent-looking'' determinant. This is also listed
as Conjecture~46 (p.~42) in Christian Krattenthaler's celebrated
essay~\cite{K} on the art of determinant evaluation.


Let, as usual, $\delta(\alpha,\beta)$ be
the Kronecker delta function
($\delta(\alpha,\beta)=1$ when $\alpha=\beta$ and
$\delta(\alpha,\beta)=0$ when $\alpha \neq \beta$), and let,
also as usual,
$$
\qbinom a b =
\frac{(1-q^a)(1-q^{a-1}) \cdots (1-q^{a-b+1})}{(1-q^b)(1-q^{b-1}) \cdots (1-q)}.
$$
Define the discrete function $a(i,j)$ by:
$$
a(i,j)=  q^{i+j-1} \left (
\qbinom{i+j-2}{i-1} + q \qbinom{i+j-1}{i} \right )+
(1+q^i)\delta(i,j)-\delta(i,j+1).
$$
Soichi Okada (\cite{O},~see also \cite{K}, Conj.~46) proved that the
$q$-TSPP conjecture is true if
$$
\det(a(i,j))_{1 \leq i,j \leq n}=
\prod_{ 1 \leq i \leq j \leq k \leq n} 
{
\left ( \frac{ 1-q^{i+j+k-1}}{1-q^{i+j+k-2}}   \right )^2
}.
$$
So in order to prove the $q$-TSPP conjecture, all we need is
to prove Okada's conjectured determinant evaluation.

\section{Certificates for Determinant Identities}

In~\cite{Z3}, an empirical (yet fully rigorous!) approach is described to
(symbolically!)
evaluate determinants $A(n):=\det (a(i,j))_{1 \leq i, j \leq n}$,
where $a(i,j)$ is a holonomic discrete function of $i$ and $j$. 
Note that this is an \emph{approach}, not a method! It is not
guaranteed to always work (and probably usually doesn't!).

Let's first describe this approach in more general terms, not just
within the holonomic ansatz.

Suppose that $a(i,j)$ is given ``explicitly'' (as it sure is here),
and we want to prove for \emph{all} $n \geq 1$ that
$$
\det (a(i,j))_{1 \leq i,j \leq n} = \Nice(n),
$$
for some \emph{explicit} expression $\Nice(n)$ (as it sure is here).

The approach is to
\emph{pull out of the hat} another ``explicit''
(possibly in a much broader sense of the word \emph{explicit})
discrete function $B(n,j)$, and then check the
identities
\begin{alignat*}3
\sum_{j=1}^n B(n,j)a(i,j)&=0, \qquad&&(1 \leq i <n < \infty), \tag{\!\textit{Soichi}\kern-.2ex}\label{eq:Soichi}\\
B(n,n)&=1, \qquad&&(1 \leq n < \infty). \tag{\!\textit{Normalization}\kern-.2ex}\label{eq:Normalization}
\end{alignat*}
If we could do that, then by \emph{uniqueness}, it would follow that
$B(n,j)$ equals the co-factor of the $(n,j)$ entry of
the $n \times n$ determinant divided by the $(n-1)\times(n-1)$ determinant
(that is the  co-factor of the $(n,n)$ entry in the $n \times n$ determinant).
Finally one has to check the identity
\begin{alignat*}1
\sum_{j=1}^{n} B(n,j) a(n,j)= \Nice(n)/\Nice(n-1) \quad(1\leq n<\infty)\tag{\!\textit{Okada}\kern-.2ex}\label{eq:Okada}
\end{alignat*}
If the suggested function $B(n,j)$ does satisfy \eqref{eq:Soichi}, \eqref{eq:Normalization}, 
and~\eqref{eq:Okada}, then the determinant identity follows as a consequence. 
So in a sense, the explicit description of $B(n,j)$ plays the role of a \emph{certificate}
for the determinant identity.

\section{The $q$-Holonomic Ansatz}

In what sense might $B(n,j)$ be explicit? In \cite{Z3} the focus was on holonomic 
sequences, in the present situation we will work with $q$-holonomic sequences.
A univariate sequence $F(n)$ is called $q$-holonomic if it satisfies a linear recurrence
of the form
\[
  a_r(q,q^n) F(n+r) + a_{r-1}(q,q^n)F(n+r-1) + \cdots + a_1(q,q^n)F(n+1) + a_0(q,q^n)F(n) = 0
\]
where $a_0,\dots,a_r$ are certain polynomials. A key feature is that $F(n)$ is 
uniquely determined by such a recurrence and the initial values $F(1),\dots,F(r)$. It is
therefore fair to accept recurrence plus initial values as an \emph{explicit description}
of the sequence~$F(n)$.

A bivariate $q$-holonomic sequence $F(n,m)$ is uniquely determined by a linear recurrence
of the form
\[
  a_r(q,q^n,q^m) F(n+r,m) + \cdots + a_1(q,q^n,q^m)F(n+1,m) + a_0(q,q^n,q^m)F(n,m) = 0
\]
where $a_0,\dots,a_r$ are certain rational functions and $F(1,m),\dots,F(r,m)$ are 
$q$-holonomic as univariate sequences in~$m$. This construction can be continued to discrete
functions of any number of variables. 

Note that while every $q$-holonomic discrete function can be described
as above, not every function that is described as above, with
\emph{arbitrary} polynomials~$a_r$ is necessarily holonomic (usually it isn't!).
However there are efficient algorithms for deciding whether a candidate
discrete function given as above is holonomic or not. 
One empirical way of doing this is to use the description
to crank out many values, and then ``guess'' a pure recurrence with polynomial
coefficients in the other variable, $m$, that can be 
routinely proved \emph{a posteriori}.

Just as holonomic sequences~\cite{Z1}, $q$-holonomic sequences have a number of important 
properties. We recall the most important ones:
\begin{enumerate}
\item If $F(n_1,\dots,n_d)$ and $G(n_1,\dots,n_d)$ are ($q$-)holonomic, 
  then so are the sequences 
  \[
    F(n_1,\dots,n_d)+G(n_1,\dots,n_d)
    \quad\text{and}\quad
    F(n_1,\dots,n_d)G(n_1,\dots,n_d).
  \]
  A ($q$-)holonomic description of these can be computed algorithmically given
  ($q$-)holonomic descriptions of $F$ and~$G$.
\item If ($q$-)holonomic descriptions of some sequences $F(n_1,\dots,n_d)$ and $G(n_1,\dots,n_d)$ 
  are given, then it can be decided algorithmically whether $F=G$.
\item If $F(n_1,\dots,n_d)$ is ($q$-)holonomic, then 
  \[
   G(n_1,\dots,n_{d-1})=\sum_{k=-\infty}^\infty F(n_1,\dots,n_{d-1},k)
  \]
  is ($q$-)holonomic.

  A ($q$-)holonomic description of $G(n_1,\dots,n_{d-1})$ can be computed algorithmically 
  given a ($q$-)holonomic descriptions of $F(n_1,\dots,n_d)$.
\end{enumerate}

\section{The Computational Challenge}

Denote by $B'(n,j)$ the sequence defined by \eqref{eq:Soichi} and \eqref{eq:Normalization}.
Why can we expect that $B'(n,j)$ is $q$-holonomic? 
\emph{A priori} there is no reason why it should be. We have to \emph{hope}.
And we can systematically \emph{search} for a potential $q$-holonomic description of~$B'(n,j)$.
If we find something, we have won, if not, we have lost,
but there is always the hope that a further search, with larger
parameters, would be successful.

One can use \eqref{eq:Soichi} and \eqref{eq:Normalization} to compute the values $B'(n,j)$ for,
say, $1\leq j\leq n\leq 35$ and then make an \emph{ansatz} for a linear recurrence, say,
of the form
\[
 \sum_{\gamma=0}^{10}
 \Bigl(
 \sum_{\beta=0}^7
 \sum_{\alpha=0}^4
 c_{\alpha,\beta,\gamma} q^{\alpha n} q^{\beta j}\Bigr) B'(n,j+\gamma) = 0.
\]
For each specific choice of $n$ and~$j$, this equation reduces to a linear equation for the 
undetermined coefficients~$c_{\alpha,\beta,\gamma}$. With different choices of $j$ and~$n$, 
we create an overdetermined linear system for the~$c_{\alpha,\beta,\gamma}$. 
Solutions of that system, if there are any, are with good probability valid recurrences 
for~$B'(n,j)$. 

So \emph{in principle,} we just have to solve a linear system. 
But \emph{in practice,} this is not as easy as it might seem.
The values of $B'(n,j)$ are rational functions in~$q$, and so are the solutions~$c_{\alpha,\beta,\gamma}$. 
A dense linear system over the rational functions with~440 unknowns cannot be solved directly with
Gaussian elimination. The intermediate expression swell would quickly blow up the matrix coefficients
to an astronomic size. Also the computation of~465 values of $B'(n,j)$ via \eqref{eq:Soichi} and 
\eqref{eq:Normalization} is not entirely for free, because it too requires solving dense linear
systems whose coefficients are rational functions in~$q$. 

We solved the system using homomorphic images. In a first step, we computed the values of $B'(n,j)$ 
with $q$ set to~$2$, and reduced modulo the prime~$p:=2^{31}-1$. This can be done quickly. 
Also the linear system for the ansatz above can be solved quickly within the finite field with
$p$ elements.
It turned out that there is a one dimensional solution space. 
At this point, there is good evidence that the $B'(n,j)$ satisfy a recurrence of the above form,
but we do not know the explicit form of the coefficients $c_{\alpha,\beta,\gamma}$ yet.
Only their homomorphic images are known. 

In the homomorphic image, 110 of the 440 coefficients $c_{\alpha,\beta,\gamma}$ are zero.
We next refined the ansatz for the recurrence by discarding the terms 
$q^{\alpha n} q^{\beta j} B'(n,j+\gamma)$ for which $c_{\alpha,\beta,\gamma}$ was found to be 
zero in the homomorphic image.
Next we repeated the computation of the $B'(n,j)$ and the solution of the linear system 
for $q=3,4,5,6,\dots,150$, always computing modulo~$p$. 
The modular images were then combined via polynomial interpolation, rational function 
reconstruction, and rational number reconstruction~\cite{vzGG} to coefficients which are rational
functions in $q$ over the integers. 

The resulting candidate recurrence has a number of remarkable features. 

\begin{enumerate}
\item The recurrence was obtained as a solution of a dense overdetermined linear system.
  An artefact solution to an overdetermined system appears only with very low probability.
\item The integer coefficients in the rational functions $c_{\alpha,\beta,\gamma}$ 
  do not exceed $43$ in absolute value. 
  For an artefact solution, integers with absolute value up to $\sqrt{p}\approx 10^9$ would 
  be expected with very high probability. 
\item The polynomials
  \[
   \sum_{\beta=0}^7 \sum_{\alpha=0}^4
    c_{\alpha,\beta,\gamma} q^{\alpha n} q^{\beta j} \qquad(\gamma=0,\dots,10)
  \]
  factorize into low degree factors. For example, the leading coefficient of
  the recurrence ($\gamma=10$) factors as
  \[
    (q^{j+6}-1) (q^{j+10}+1) (q^n-q^{j+9}) (q^n-q^{j+10}) (q^{j+n+9}-1) (q^{j+n+10}-1).
  \]
  For an artefact solution, it would be expected with very high
  probability that all the polynomials are irreducible.
\item The recurrence produces the correct terms of $B'(n,j)$ for values $35<n\leq 200$ 
  for $q=2$ and modulo~$p$, although these terms were not used in the computation 
  of the recurrence. 
  For an artefact solution, this is expected to happen with very low probability only.
\item The recurrence produces the correct terms of $B'(n,j)$ for small $n$ and~$j$
  if $q$ is left symbolic or set to a numeric value different from $2,3,\dots,100$,
  For an artefact solution, this is expected to happen with very low probability only.
\end{enumerate}

We have not the slightest doubt that the recurrence we found is correct. 
For a rigorous proof, we can define (``pull out of the head'') a sequence $B(n,j)$ 
by a $q$-holonomic description consisting of the recurrence we discovered and some suitable 
univariate recurrences and initial values that are easy to obtain.
It is contained in the Maple package {\tt qTSPP} accompanying this article. 
The much easier $q=1$ case (that would give a new proof to the already proved Stembridge theorem)
is contained in the Maple package {\tt TSPP}.
Proving that $B(n,j)=B'(n,j)$
amounts to proving that $B(n,j)$ satisfies \eqref{eq:Soichi}. (Equation \eqref{eq:Normalization}
is automatically satisfied.) Thanks to algorithms of Chyzak, Salvy, Takayama~\cite{CS,T}, proving 
\eqref{eq:Soichi} \emph{in principle} reduces to finitely many routine calculations.

Finally, \eqref{eq:Okada} is also of the form $A=B$ where both
sides are $q$-holonomic. The left side is $q$-holonomic because
of the closure under multiplication and definite-summation,
and the right side, $\Nice(n)/\Nice(n-1)$ is not just $q$-holonomic
(a solution of \emph{some} linear recurrence  with polynomial coefficients
(in $q,q^n$)\,) but in fact \emph{closed-form}
(the defining recurrence is first-order).

Summarizing, we have found a certificate~$B(n,j)$ for Okada's conjectured determinant 
identity
\[
 \det(a(i,j))=\prod_{1\leq i\leq j\leq k\leq n}\left(\frac{1-q^{i+j+k-1}}{1-q^{i+j+k-2}}\right)^2.
\]
It only remains to prove rigorously that our certificate really is a certificate. 
While such a proof can \emph{in principle} be carried out automatically, we found that
\emph{in practice,} i.e., with the currently available algorithms, software, and hardware,
it remains a computational challenge. Even for the case $q=1$, we are at present unable 
to complete the necessary non-commutative Gr\"obner basis computations.

\section{The Semi-Rigorous Shortcut}

We believe that even today, performing the computations for a rigorous proof is feasible,
but it would require a huge technical effort. But why bother?
First, if we wait for twenty more years, the availability of better algorithms, 
better software, and better hardware will probably
enable us to finish up these finitely many routine calculations
with no sweat. Besides, since now we know for sure that
a \emph{fully rigorous} proof exists, do we really want to
see it? It won't give us any new \emph{insight}. The beauty of
the present approach is in the \emph{meta-insight}, reducing 
the statement of the conjecture
to a \emph{finite} calculation. Furthermore, we know \emph{a priori} that
there exists an operator $P(q,q^n,N)$ (where $N$ is the shift operator in $n$:
$Nf(n):=f(n+1)$) that annihilates the
difference of the left and right sides of \eqref{eq:Okada}. If that
operator has order $L$, say, then a completely rigorous proof
would be to check \eqref{eq:Okada} for $1 \leq n \leq L$.
At present, we are unable to find $P$, and hence do not know
the value of $L$. But it is very reasonable that $L$ would be
less than, say, $400$, and checking the first $400$ cases
of \eqref{eq:Okada} (and analogously for \eqref{eq:Soichi}) is certainly
doable (we did it for $L=100$, and $L=400$ for TSPP, but
you are welcome to go further). These are done in procedures
{\tt CheckqTSPP} in the Maple package {\tt qTSPP}, and
{\tt CheckTSPP} in the Maple package {\tt TSPP}, respectively.
The corresponding input and output files can be found
in the webpage of this article mentioned above.
As a technical aside, let's confess that Maple running on our computer
was only able to check \eqref{eq:Soichi} and \eqref{eq:Okada} for $L \leq 30$, for
\emph{symbolic}~$q$, but for random numerical choices of $q$ it went
up to $L=100$, and it is easy to see that with sufficiently many
choices of numerical~$q$ for a given $L$, one can prove it for
symbolic~$q$.

In 1993, Zeilberger~\cite{Z2} proposed the notion of
\emph{semi-rigorous} proof. At the time he didn't have any
\emph{natural} examples. The present determinant evaluation,
that was shown by Okada to imply a long-standing open
problem in enumerative combinatorics, is an \emph{excellent}
example of a semi-rigorous proof that is (at least) as good 
as a rigorous proof. Let us conclude
by promising that if any one is willing
to pay us \$$10^7$ 
(ten million US dollars), we will be more than glad to fill in the details.

\section{Postscript}

This article was originally submitted to a special volume 
of the {\it Seminaire Lotharingien de Combinatoire} (SLC)
in memory of Pierre Leroux. While the editors and referees were
willing to accept our paper, they demanded that we change the
title and ``tone down'' our claim that we  have a proof
(even modulo a finite amount of calculations). Since there is
a ``mathematical'' possibility (as the French would put it) that
our ``proof plan'' would not work out, in which case we {\it have nothing}.

We agree that there is a positive probability that our proof would turn out
to be wrong. But that probability is orders-of-magnitude smaller
than the probability that the editors of SLC do not
exist. After all they, along with all of us, may be characters in a dream of
a giant, and we would all disappear once that giant wakes up.

Since we believe that our title is a good one, and all our claims are sound,
we decided to withdraw the paper from SLC, and publish it in the
much more enlightened journal {\it The Personal Journal of Ekhad and Zeilberger},
as well as in {\it arxiv.org}.

However, as a concession to the human sentiments expressed by the editors and referees
and for those impatient people who can't wait twenty years, or cannot
afford ten million US dollars, let us conclude with a sketch on
how to hopefully make the present approach yield a fully rigorous proof
with today's software and hardware.

The key is to take advantage of the special structure of the
entries $a(i,j)$, and not just consider them as yet-another-holonomic
sequence. For the sake of simplicity, let's consider the $q=1$ case.
Analogous considerations apply to the $q$-case.

When $q=1$, the matrix entry is:
$$
a(i,j)=  
\binom{i+j-2}{i-1} + \binom{i+j-1}{i} +
2\delta(i,j)-\delta(i,j+1).
$$
Let's write it as:
$$
a(i,j)=  a'(i,j)+
2\delta(i,j)-\delta(i,j+1) \quad ,
$$
where
$$
a'(i,j)=  
\binom{i+j-2}{i-1} + \binom{i+j-1}{i} \quad .
$$
It is also helpful to define $C(n,j):=B(n,n-j)$.
Note that $C(n,j)$ is defined to be $0$ for $j<0$,
$C(n,0)=1$ and $C(n,1)$ has a certain conjectured
holonomic description as a sequence in $n$.

At this point it may be fruitful to introduce the
sequence of polynomials
$$
f_n(x)=\sum_{j=0}^{n} C(n,j)x^j \quad .
$$
(it may be more efficient to let the sum range from $j=0$ to
$j=n+2$).
The $j$-free recurrence for $B(n,j)$ given in the
package {\tt TSPP} translates to a certain
linear recurrence equation, in $n$,
with polynomial coefficients in $(n,x)$, for $f_n(x)$,
and the original $N$-free recurrence, translates to a certain
{\it linear differential equation}, in $x$, with  polynomial
coefficients in $n$ and $x$. Just like all the 
classical orthogonal polynomials (Legendre, Laguerre, Hermite, Jacobi, etc.),
except that the relevant equations are no longer second-order,
(and the $\{f_n(x)\}$ are not orthogonal).
In particular the discrete-continuous function $(n,x) \rightarrow f_n(x)$
has a full holonomic description in its arguments.

Now both $(Soichi)$ and $(Okada)$ can be easily transcribed to certain
simply-stated constant-term identities in $(n,x)$.

Indeed, $(Soichi)$ is equivalent to the constant-term identity
$$
CT \left [
\left \{ {{x(2-x)} \over {(1-x)^{i+1}} } + 2x^i-x^{i-1} \right \} {{f_n(x)} \over {x^n}} \right ] \, = \, 0 \quad,
\eqno(Soichi')
$$
where $CT$ stands for ``constant term'', i.e. ``coefficient of $x^{0}$''. Calling the 
left side $L(n,i)$, note that we can induct on {\it both} $n$ and $i$, so we would be done
once we have an annihilating operator of the form $P(N,I,n,i)$ (here $N$ is shift in the
discrete variable $n$, and $I$ is shift in the discrete variable $i$).  
Calling the {\it constant-termand} $F(n,i,x)$, this means that all we
need is find an operator of the form $P(N,I,n,i,x{{d} \over {dx}})$
annihilating $F(n,i,x)$. It shouldn't be too hard to get two operators
annihilating $F(n,i,x)$ out of the ones that we already have for $f_n(x)$, and then 
express them as $P(n,i,x,N,I,x{{d} \over {dx}})$, $Q(n,i,x,N,I,x{{d} \over {dx}})$, say,
and use non-commutative Gr\"obner base eliminiation to eliminate $x$.
This may be much more efficient than finding ``pure'' operators of the form $P(N,n,i,x{{d} \over {dx}})$,
or $P(I,n,i,x{{d} \over {dx}})$,
where we have to eliminate two ``variables'' in the non-commutative
algebra of recurrence-differential operators $C[n,i,x,N,I,x{{d} \over {dx}}]$.

As for $(Okada)$, it is equivalent to
$$
CT \left [
\left \{ {{x(2-x)} \over {(1-x)^{n+1}} } + 2x^n-x^{n-1} \right \} {{f_n(x)} \over {x^n}} \right ] \, = \, \Nice(n)/\Nice(n-1) \quad .
\eqno(Okada')
$$

It is very possible, that going via this {\it continuous-discrete} route would make
the problem tractable with today's software and hardware,
and we leave it as a challenge in case any of our readers is
interested enough to convert our (currently) semi-rigorous proof into
a fully rigorous proof, rather than  wait patiently for twenty years.

\end{document}